\RequirePackage{fix-cm}
\documentclass{svjour3}                     
\smartqed  
\usepackage{graphicx}
\usepackage[numbers,sort&compress]{natbib}
 \usepackage{mathptmx}      
\usepackage{latexsym}
\begin{document}

\title{Immersions into Statistical Manifolds
}


\author{Mahesh T V         \and K S Subrahamanian Moosath$^{*}$\thanks{ *corresponding author}
        }


\institute{Department of Mathematics \at
              Indian Institute of Space Science and Technology \\
              Thiruvanathapuram, Kerala, India-695547\\
              \email{maheshtv.16@res.iist.ac.in(Mahesh T V)} \\ 
              \email{smoosath@iist.ac.in (K S Subrahamanian Moosath)}         
           }

\date{Received: date / Accepted: date}

\maketitle

\begin{abstract}
This paper studies the geometry of immersions into statistical manifolds. A necessary and sufficient condition is obtained for statistical manifold structures to be dual to each other for a non-degenerate equiaffine immersion. Then we obtain conditions for realizing an $n$ dimensional statistical manifold in an $(n+1)$-dimensional statistical manifold and its converse. Centro-affine immersion of codimension two into a dually flat statistical manifold is defined. Also we have shown that statistical manifold realized in a dually flat statistical manifold of codimension two is conformally-projectively flat.
\keywords{Statistical Manifold\and Statistical Immersion}
 \subclass{MSC 53A15 \and MSC 53C42}
\end{abstract}
\section{Relevance Of The Work}
Information geometry is the result of applying non-Euclidean geometry to probability theory in general and statistical inference and estimation in particular. It has got applications in information science, artificial neural networks, signal processing, statistical mechanics, neuroscience etc. One of the possibilities of deforming the geometric structures of the statistical manifolds is using the affine immersion. In this paper we discuss the geometry of immersions into statistical manifolds. Mainly the conditions for realizing $n$- dimensional statistical manifolds in $(n+1)$- dimensional statistical manifolds and in dually flat statistical manifolds of dimension $(n+2)$.
\section{Introduction}
\label{intro}
Statistical manifolds have close relation to the geometry of immersions.
  Amari posed the problem of finding conditions for realization of  statistical manifolds  in  affine spaces. Kurose \cite{kurose1990dual} proved a necessary and sufficient condition for a statistical manifold to be realized by affine immersions of codimension one. Matsuzoe \cite{matsuzoe1998realization} obtained a necessary and sufficient condition for a statistical manifold to be realized by a non-degenerate centro-affine, equiaffine immersion into an affine space of codimension two. In this paper, we consider immersions into statistical manifolds as well as statistical immersions and study the geometry of such submanifolds.  \par
  
 In section 2, first we give basic ideas on statistical manifolds and immersions into statistical manifolds. Apart from the induced statistical structure $(f^{*}\tilde{g}, \nabla,\nabla^{*} )$ we have the inherited statistical structures $(h,\nabla)$ and $(h^{*},\nabla^{*})$ for a non-degenerate equiaffine immersion. The induced connections $\nabla$, $\nabla^{*}$ are dual to each other with respect to the induced metric $f^{*}\tilde{g}$. We prove a necessary and sufficient condition for the inherited statistical manifold structures to be dual to each other.

 In section 3, statistical immersions and centro-affine immersions of codimension two into dually flat statistical manifolds are considered. Kurose \cite{kurose1990dual} gave a necessary and sufficient condition for realizing a statistical manifold by an affine immersion of codimension one. We have given conditions for realizing a statistical manifold of dimension $n$ in a statistical manifold of dimension $(n+1)$. In the case of statistical immersion we prove its converse for statistical manifolds of dimension $n\geq 3$. Then we define centro-affine immersion of codimension two into a dually flat statistical manifold. Matsuzoe \cite{matsuzoe1998realization} showed that a statistical manifold of dimension $n$ realized in $\mathbf{R}^{n+2}$ is conformally-projectively flat and also obtained its converse. For a non-degenerate centro-affine, equiaffine immersion of a manifold of dimension $n$ into a dually flat statistical manifold of dimension $(n+2)$ we prove that the inherited statistical manifold is conformally- projectively flat.

  Throughout this paper we assume that all the objects are smooth and affine connections are torsion-free.
\section{ Immersions into Statistical manifolds}
 
In this section first we give a brief description on statistical manifolds \cite{amari2016information}, \cite{amari2000method} and immersions into statistical manifolds. Then we prove a necessary and sufficient condition for the inherited statistical manifold structures to be dual to each other. 
 
 A pseudo - Riemannian manifold $(\mathbf{M},g)$ with a torsion free affine connection $\nabla$ is called a statistical manifold if $\nabla g$ is symmetric. For a statistical manifold $(\mathbf{M},\nabla,g)$ the dual connection  $\nabla^{*}$ is defined by  
  \begin{equation}
   Xg(Y,Z) = g(\nabla_{X}Y,Z) + g(Y,\nabla^{*}_{X} Z)
  \end{equation}
 for  vector fields $X,Y$ and $Z$ in $\mathcal{X}(\mathbf{M})$, where $\mathcal{X}(\mathbf{M})$ denotes the set of all vector fields on $\mathbf{M}$. If $(\nabla,g)$ is a statistical structure on 
 $\mathbf{M}$ so is $(\nabla^{*},g)$. Then $(\mathbf{M},\nabla^{*},g)$ becomes a statistical manifold  called the dual statistical manifold of $(\mathbf{M},\nabla,g)$.  Let $R^{\nabla}$ and $R^{\nabla^{*}}$ be the curvature tensors of $\nabla$ and $\nabla^{*}$, respectively. It follows from $(1)$ that 
 \begin{equation}
 g(R^{\nabla}(X,Y)Z,W) = -g(Z,R^{\nabla^{*}}(X,Y)W)
 \end{equation}
for  $X,Y,Z$ and $W$ in  $\mathcal{X}(\mathbf{M})$. We say $(\mathbf{M},\nabla,\nabla^{*},g)$ has constant curvature $k$ if
 \begin{equation}
 R^{\nabla}(X,Y)Z  = k \lbrace g(Y,Z)X - g(X,Z)Y \rbrace.
 \end{equation}
 A statistical manifold with  curvature zero is called a flat statistical manifold and in that case $(\mathbf{M},\nabla,\nabla^{*},g)$  is called a dually flat statistical manifold.\par
 
 Two statistical manifolds $(\mathbf{M},\nabla,g)$ and $(\mathbf{M},\tilde{\nabla},\tilde{g})$ are said to be conformally-projectively equivalent if there exist two positive functions $\phi$ and $\psi$ on $\mathbf{M}$ such that 
\begin{eqnarray*}
\tilde{g}(X,Y)&=&\phi \psi g(X,Y),\\
g(\tilde{\nabla}_{X}Y ,Z)&=&g(\nabla_{X}Y,Z)-d(log\phi)(Z)g(X,Y)\\
&+& d(log\psi)(X)g(Y,Z)+ d(log\psi)(Y)g(X,Z).
\end{eqnarray*}
A statistical manifold $(\mathbf{M},\nabla,g)$  is said to be conformally-projectively flat if it is conformally-projectively equivalent to a flat statistical manifold in a neighbourhood of an arbitrary point of $\mathbf{M}$.\\ 

 Let $\mathbf{M}$ be an  $n$-dimensional manifold and $(\tilde{\mathbf{M}},\tilde{\nabla},\tilde{g})$ be an   $(n+1)$-dimensional statistical manifold. Let $f:\mathbf{M}\longrightarrow \tilde{\mathbf{M}}$ be an immersion. The induced metric $f^{*} \tilde{g}$ and the induced connection $\nabla$ on $\mathbf{M}$ are defined  as follows
 \begin{eqnarray}
 f^{*} \tilde{g}(\nabla_{X}Y, Z) &=& \tilde{g}(\tilde{\nabla}_{X}{f_{*}Y},f_{*}Z)
\end{eqnarray}
for any $X,Y,Z \in \mathcal{X}(\mathbf{M})$.

Note that
\begin{eqnarray*}
(\nabla_{X}f^{*} \tilde{g})(Y,Z) &=& Xf^{*} \tilde{g}(Y,Z) - f^{*} \tilde{g}(\nabla_{X}Y,Z)- f^{*} \tilde{g}(Y,\nabla_{X}Z)\\
                   &=& X\tilde{g}(f_{*}Y,f_{*}Z)-\tilde{g}(\tilde{\nabla}_{f_{*}X}f_{*}Y,f_{*}Z)-\tilde{g}(\tilde{f_{*}Y,\nabla}_{f_{*}X}f_{*}Z)\\
                   &=& (\tilde{\nabla}_{X} \tilde{g})(f_{*}Y,f_{*}Z).
\end{eqnarray*}
Thus $(\mathbf{M},\nabla,f^{*} \tilde{g})$ is also a statistical manifold.\par

 Let $f:\mathbf{M}\longrightarrow \tilde{\mathbf{M}}$ be an immersion of codimension one with unit normal vector field $\xi$ along $f$. Then 
for each $p \in \mathbf{M} $
\begin{equation}
T_{f(p)}(\mathbf{\tilde{\mathbf{M}}}) = f_{*}(T_{p}(\mathbf{M})) + span\lbrace\xi_{p}\rbrace.
\end{equation}
Also the Gauss and Weingarten formulae \cite{nomizu1987geometry}, \cite{vos1989fundamental} are
\begin{enumerate}
\item $\tilde{\nabla}_{X}f_{*}Y$ = $f_{*}(\nabla^{\top}_{X}Y)+h(X,Y)\xi$
\item $\tilde{\nabla}_{X}\xi$ = $-f_{*}(SX)+\tau(X)\xi$
\item $\tilde{\nabla}^{*}_{X}f_{*}Y$ = $f_{*}(\nabla^{\top*}_{X}Y)+h^{*}(X,Y)\xi$
\item $\tilde{\nabla}^{*}_{X}\xi$ = $-f_{*}(S^{*}X)+\tau^{*}(X)\xi$
\end{enumerate}
for $X,Y \in \mathcal{X}(\mathbf{M})$, where $\tilde{\nabla}^{*}$ is the dual connection of  
$\tilde{\nabla}$ with respect to $\tilde{g}$, $h(X,Y)$ and $h^{*}(X,Y)$ are  symmetric bilinear forms on tangent space $T_{p}(\mathbf{M})$ for $p$ in $\mathbf{M}$.  $S$ and $S^{*}$ are tensor fields of type $(1,1)$ and $\tau,\tau^{*}$ are $1$-forms. We call $S$ ($S^{*}$) the shape operator and  $\tau$ ($\tau^{*}$) the transversal connection form for $f$ and  the induced connections $\nabla^{\top}$ ($\nabla^{\top*}$).\\ 

The  connection $\nabla^{\top}$ induced by Gauss formula  coincides with the induced connection $\nabla $ on $\mathbf{M}$. So we write $\nabla $ instead of $\nabla^{\top}$. Thus on $\mathbf{M}$ we have induced connections $\nabla$ and $\nabla^{*}$ and affine fundamental forms $h$ and $h^{*}$. Using the Gauss and Weingarten formulae we have the following observations:\\

\noindent 1. Let $\mathbf{M}$ be an  $n$-dimensional manifold, $(\tilde{\mathbf{M}},\tilde{\nabla},\tilde{g})$ be an $(n+1)$-dimensional statistical manifold and let $f:\mathbf{M}\longrightarrow \tilde{\mathbf{M}}$ be an immersion of codimension one. Then the induced connections $\nabla$ and $\nabla^{*}$ on $\mathbf{M}$ are dual with respect to the induced metric $f^{*}\tilde{g}$.\\

\noindent 2. Let $f:\mathbf{M}\longrightarrow \tilde{\mathbf{M}}$ be an immersion of codimension one. Then for each $X,Y \in \mathcal{X}(\mathbf{M})$
\begin{enumerate}
  \item $h(X,Y)$ = $f^{*} \tilde{g}(S^{*}X,Y)$
  \item $h^{*}(X,Y)$ = $f^{*} \tilde{g}(SX,Y)$ 
  \item $\tau(X)+\tau^{*}(X) $ = $0$
  \end{enumerate}
Now suppose $(\tilde{\mathbf{M}},\tilde{\nabla},\tilde{g})$ has constant curvature $\tilde{k}$, then the fundamental equations \cite{nomizu1987geometry}, \cite{vos1989fundamental} are
 \begin{eqnarray*}
 R^{\nabla}(X,Y)Z &=& \tilde{k}\lbrace f^{*} \tilde{g}(Y,Z)X - f^{*} \tilde{g}(X,Z)Y \rbrace +  h(Y,Z)SX\\&& - h(X,Z)SY ~~~~~~~~~~~~~(Gauss) \\
(\nabla_{X}h )(Y,Z)+\tau(X)h(Y,Z)&=&(\nabla_{Y}h )(X,Z)+\tau(Y)h(X,Z)~~~(Codazzi \hspace{.1cm}
for \hspace{.1cm} h) \\
  (\nabla_{X}S)(Y)-\tau(X)SY &=&  (\nabla_{Y}S)(X)-\tau(Y)SX ~~~~~~~~~~(Codazzi \hspace{.1cm} for \hspace{.1cm} S)\\
 h(X,SY)-h(SX,Y) &=& d\tau(X,Y) ~~~~~~~~~~~~~~~(Ricci). 
 \end{eqnarray*}
 where $R^{\nabla}(X,Y)Z$ denotes the curvature tensor with respect to $\nabla$ in $\mathbf{M}$. Similarly one can write the equations with respect to dual connection also.
 
 \begin{definition}
 Let $f:\mathbf{M}\longrightarrow \tilde{\mathbf{M}}$ be an immersion of codimension one. Then $f$ is said to be non-degenerate if $h$ is non-degenerate and $f$ is  equiaffine if $\tau$ = $0$. 
 \end{definition}
 \begin{remark}
  For a non-degenerate equiaffine immersion $f:\mathbf{M}\longrightarrow \tilde{\mathbf{M}}$ with $\tilde{\mathbf{M}}$ has constant curvature, from the Codazzi equations for $h$ and $h^{*}$, it is obvious that both $(\mathbf{M},\nabla,h)$ and $(\mathbf{M}, \nabla^{*}, h^{*})$ become statistical manifolds.
 \end{remark}
 
  \begin{definition}
  Statistical manifolds $(\mathbf{M}, \nabla, h)$ and $(\mathbf{M}, \nabla^{*}, h^{*})$ are said to be dual to each other if $h$ = $h^{*}$ and $\nabla$, $\nabla^{*}$ are dual with respect to $h$.
 \end{definition}
 
  Now we prove a necessary  and sufficient condition for the inherited statistical manifolds $(\mathbf{M},\nabla, h)$ and $(\mathbf{M},\nabla^{*}, h^{*})$ to be dual to each other.
\begin{theorem}
Let  $\mathbf{M}$ be an $n$-dimensional manifold and $(\tilde{\mathbf{M}},\tilde{\nabla},\tilde{g})$ be an $(n+1)$-dimensional statistical manifold with constant curvature $\tilde{k}$. Let  $f:\mathbf{M}\longrightarrow \tilde{\mathbf{M}}$ be a non-degenerate, equiaffine  immersion of codimension one. Then  $(\mathbf{M},\nabla,h)$ and  $(\mathbf{M},\nabla^{*},h^{*})$ are dual to each other if and only if $S$ = $S^{*}$ = $\lambda I$ for some constant $\lambda$. Moreover $h= \lambda f^{*}\tilde{g}$.
\end{theorem}
\begin{proof}
Suppose  $S$ = $S^{*}$ = $\lambda I$ for some constant $\lambda$. From the above observation$(2)$
 \begin{eqnarray*}
 h(X,Y) &=& f^{*}\tilde{g}(S^{*}X,Y)\\
        &=& \lambda f^{*}\tilde{g}(X,Y)
 \end{eqnarray*}
 Similarly $h^{*}(X,Y)$ = $\lambda f^{*}\tilde{g}(X,Y)$, so we have $h$ = $h^{*}$ = $\lambda f^{*}\tilde{g}$. Since $\nabla$ and $\nabla^{*}$ are dual with respect to $f^{*}\tilde{g}$, the statistical manifolds $(\mathbf{M},\nabla,h)$ and  $(\mathbf{M},\nabla^{*},h^{*})$ are dual to each other.\par Conversely, let $(\mathbf{M},\nabla,h)$ and  $(\mathbf{M},\nabla^{*},h^{*})$ be dual to each other, then $h$ = $h^{*}$ , $\nabla$ and $\nabla^{*}$ are dual with respect to $h$. So 
 \begin{equation}
 Zh(X,Y) = h(\nabla_{Z}X , Y) + h(X, \nabla^{*}_{Z}Y).
 \end{equation}
  Now consider
\begin{eqnarray}
  Zf^{*}\tilde{g}(SX,Y) &=& f^{*}\tilde{g}(\nabla_{Z}SX , Y) + f^{*}\tilde{g}(SX, \nabla^{*}_{Z}Y) \nonumber \\
          &=& f^{*}\tilde{g}(\nabla_{Z}SX , Y) + h(X, \nabla^{*}_{Z}Y).
\end{eqnarray}
 Since $h(X,Y)$ = $f^{*}\tilde{g}(S^{*}X,Y)$ = $f^{*}\tilde{g}(SX,Y)$ = $h^{*}(X,Y)$, from $(6)$ and $(7)$ we get
 \[ \nabla_{Z}SX = S(\nabla_{Z}X),\]
 which implies  $(\nabla_{Z}S)X$ = $0$, then $S$  = $\lambda I$ for some constant $\lambda$.\\ Therefore  $S$ = $S^{*}$ = $\lambda I$. Also note that the induced metric $h$ = $\lambda f^{*}\tilde{g}.$ 
\end{proof}

\begin{remark}
Let $(\tilde{\mathbf{M}},\tilde{\nabla},\tilde{g})$ be an   $(n+1)$-dimensional statistical manifold. If $f:\mathbf{M}\longrightarrow \tilde{\mathbf{M}}$ is an immersion of codimension one then $(\mathbf{M},f^{*} \tilde{g},\nabla,\nabla^{*})$ is a statistical manifold and if $f:\mathbf{M}\longrightarrow \tilde{\mathbf{M}}$ is a non-degenerate equiaffine immersion of codimension one then $(\mathbf{M},h,\nabla,\nabla^{*})$ is a statistical manifold if and only if $S$ = $S^{*}$ = $\lambda I$ for some constant $\lambda$.
\end{remark}

\section{Statistical Immersions}
Kurose \cite{kurose1990dual} showed that a statistical manifold of dimension $n$ can be realized in $\mathbf{R}^{n+1}$ if and only if it has constant curvature. In this section, we obtain conditions for realizing a statistical manifold of dimension $n$ in a statistical manifold of dimension $n+1$ and its converse. Then we define centro-affine immersion of codimension two  into dually flat statistical manifolds. Matsuzoe \cite{matsuzoe1998realization} showed that a statistical manifold of dimension $n$ can be realized in $\mathbf{R}^{n+2}$ if and only if it is conformally-projectively flat. We prove that a statistical manifold realized in a dually flat statistical manifold of codimension two is conformally-projectively flat.

\begin{definition}
Let $(\tilde{\mathbf{M}},\tilde{\nabla},\tilde{g})$  and  $(\mathbf{M}, \nabla, g)$ be statistical manifolds. An immersion $f:(\mathbf{M},\nabla, g)$\\ $\longrightarrow (\tilde{\mathbf{M}},\tilde{\nabla},\tilde{g}) $ is called a statistical immersion if
 \begin{eqnarray*}
                g  &=& f^{*} \tilde{g} \\
 g(\nabla_{X}Y, Z) &=& \tilde{g}(\tilde{\nabla}_{X}{f_{*}Y},f_{*}Z)
\end{eqnarray*} for any $X,Y,Z \in \mathcal{X}(\mathbf{M}).$\\

\noindent An $n$-dimensional statistical manifold $(\mathbf{M},\nabla,\nabla^{*},g)$ is said to be realized in an $(n+1)$-dimensional statistical manifold $(\tilde{\mathbf{M}},\tilde{\nabla},\tilde{\nabla}^{*},\tilde{g})$   if there exist a non-degenerate equiaffine immersion $f:\mathbf{M}\longrightarrow \tilde{\mathbf{M}}$ such that the affine fundamental form $h$ equal to $g$ and the induced connections coincide with $\nabla$ and $\nabla^{*}$.
\end{definition}
\begin{note}
Let $(\mathbf{M},\nabla,\nabla^{*},g)$ be an $n$-dimensional statistical manifold realized in an $(n+1)$-dimensional statistical manifold $(\tilde{\mathbf{M}},\tilde{\nabla},\tilde{\nabla}^{*},\tilde{g})$. Then $(\mathbf{M},\nabla,\nabla^{*},h)$ becomes a statistical manifold, where $h$ is the affine fundamental form.
\end{note}

 In the case of affine immersions of codimension one Kurose \cite{kurose1990dual}  has proved  that (i) if simply connected and connected statistical manifold of dimension $n$ has constant curvature then it can be realized in $\mathbf{R}^{n+1}$ and (ii) if a simply connected and connected statistical manifold of dimension $n\geq 3$ is realized in $\mathbf{R}^{n+1}$ then it has constant curvature. In the case of immersions into statistical manifolds of codimension one we can prove that
\begin{theorem}
If $(\mathbf{M},\nabla,\nabla^{*},g)$ is an $n$-dimensional and $(\tilde{\mathbf{M}},\tilde{\nabla},\tilde{\nabla}^{*},\tilde{g})$ is an $(n+1)$- dimensional simply connected and connected statistical manifolds with constant curvatures $k$ and $\tilde{k}$ respectively. Then $M$ is realized in $\tilde{\mathbf{M}}$.
\end{theorem}
Also we have the converse of the above result.
\begin{theorem}
Let $(\tilde{\mathbf{M}},\tilde{\nabla},\tilde{\nabla}^{*},\tilde{g})$ be an $(n+1)$-dimensional connected statistical manifold with constant curvature $\tilde{k}$ and $(\mathbf{M},\nabla,\nabla^{*},g)$  be a connected statistical manifold of dimension $n \geq 3$. If there exist a statistical immersion $(f,\xi):(\mathbf{M},\nabla,\nabla^{*},g) \longrightarrow (\tilde{\mathbf{M}},\tilde{\nabla},\tilde{\nabla}^{*},\tilde{g})$ such that $h = h^{*} = g$ then $(\mathbf{M},\nabla,\nabla^{*},g)$ has constant curvature.
\end{theorem}
\begin{proof}
Since, by Gauss equation
\begin{eqnarray*}
R^{\nabla}(X,Y)Z &=& \tilde{k}\lbrace g(Y,Z)X - g(X,Z)Y \rbrace +  h(Y,Z)SX\\&& - h(X,Z)SY
\end{eqnarray*}
it is enough to show that $S = \lambda I$ for some constant $\lambda$. Now consider
\begin{eqnarray*}
R^{\nabla^{*}}(X,Y)W &=& \tilde{k}\lbrace g(Y,W)X - g(X,W)Y \rbrace +  h^{*}(Y,W)S^{*}X\\&& - h^{*}(X,W)S^{*}Y.
\end{eqnarray*}
 $\nabla$ and $\nabla^{*}$ are dual with respect to $g$, so we have
\begin{equation}
 g(R^{\nabla}(X,Y)Z,W) = -g(Z,R^{\nabla^{*}}(X,Y)W)
 \end{equation}
for  $X,Y,Z$ and $W$ in  $\mathcal{X}(\mathbf{M})$. Then 
\begin{eqnarray}
  g(Y,Z)g(SX,W) - g(X,Z)g(SY,W) &=&  g(Y,W)g(S^{*}X,Z)\nonumber \\ &-& g(X,W)g(S^{*}Y,Z)
\end{eqnarray}
Set $L = \frac{1}{n} tr(S)$ and $L^{*} = \frac{1}{n} tr(S^{*})$. Now taking trace in $X$- and $W$- components in $(9)$ we get 
\begin{equation}
nLg(Y,Z)-g(S(Y),Z)+g(S^{*}Y,Z)-ng(S^{*}Y,Z) = 0.
\end{equation}
Again, taking trace in $Y$- and $Z$-components, we get
\begin{equation}
L=L^{*}
\end{equation}
Also from $(10)$ we have
\begin{equation}
nLI = S + (n-1)S^{*}
\end{equation}
Since equation $(9)$ is symmetric in $S$ and $S^{*}$ 
\begin{equation}
nLI = S^{*} + (n-1)S
\end{equation}
$(12)$ and $(13)$ imply $S= LI$ for $n \geq 3$. Hence $(\mathbf{M},\nabla,\nabla^{*},g)$ has constant curvature.
\end{proof}
\begin{note}
Matsuzoe\cite{matsuzoe1998realization} considered centro-affine immersions of statistical manifolds of dimension $n$ into the affine space $\mathbf{R}^{n+2}$.
\end{note}
 
Now we define centro-affine immersion into a dually flat statistical manifold of codimension two.
\begin{definition}
Let $(\tilde{\mathbf{M}},\tilde{\nabla},\tilde{\nabla}^{*},\tilde{g})$  be a dually flat statistical manifold of dimension $(n+2)$ and $\mathbf{M}$ be an $n$-dimensional manifold. Let 
   $\gamma$ = $\sum_{i=1}^{n+2}\theta^{i}\frac{\partial}{\partial \theta_{i}}$  be the radial vector field of $\tilde{\mathbf{M}}$ with respect to the affine co-ordinate $[\theta^{i}]$ of $\tilde{\nabla}$.    $f :\mathbf{M}\longrightarrow\tilde{\mathbf{M}} $ is called a centro-affine immersion of codimension two if there exist a unit normal vector field $\xi$ such that
\begin{eqnarray*}
T_{f(p)}(\tilde{\mathbf{M}}) = f_{*}(T_{p}(\mathbf{M})) + span\lbrace \xi \rbrace + span\lbrace \gamma\rbrace.
\end{eqnarray*}
\end{definition} 
 In this case the Gauss and Weingarten formulae are 
 
 \begin{eqnarray}
 \tilde{\nabla}_{X}f_{*}(Y) &=&  f_{*}(\nabla_{X}Y)+h(X,Y)\xi + T(X,Y)\gamma \\
 \tilde{\nabla}_{X}\xi      &=&  -f_{*}(SX)+\tau(X)\xi + \sigma(X)\gamma\\
 \tilde{\nabla}_{X}\gamma     &=&  -f_{*}(X)
\end{eqnarray}
where $S$ is the  affine shape operator, $\tau(X)$ and $\sigma(X)$ are transversal connection forms and 
\begin{center}
$h,T: T_{p}(\mathbf{M}) \times T_{p}(\mathbf{M}) \longrightarrow \mathbf{R}$
\end{center}
are affine fundamental forms. The Gauss and Weingarten formulae for dual connection can also be written similarly.\\

 Let $f:\mathbf{M}\longrightarrow \tilde{\mathbf{M}}$ be a centro-affine immersion of codimension two. Then $f$ is called non-degenerate if $h$ is non-degenerate and $f$ is  equiaffine if $\tau$ = $0$.

\begin{remark}
As in the case of codimension one, we can easily see form the fundamental equations of centro-affine immersion of codimension two that both $(\mathbf{M},\nabla,h)$ and $(\mathbf{M},\nabla^{*},h^{*})$ are statistical manifolds for non-degenerate equiaffine, centro affine immersion.
\end{remark}

In this case also we can show that

\begin{theorem}
Let $(\tilde{\mathbf{M}},\tilde{\nabla},\tilde{\nabla}^{*},\tilde{g})$  be a dually flat statistical manifold of dimension $(n+2)$ and $\mathbf{M}$ be a manifold of dimension $n$.
If $f: \mathbf{M}\longrightarrow \tilde{\mathbf{M}}$ is a non-degenerate centro-affine, equiaffine immersion of codimension two then $(\mathbf{M},\nabla,h)$ and  $(\mathbf{M},\nabla^{*},h^{*})$ are dual to each other if and only if $S$ = $S^{*}$ = $\lambda I$ for some constant $\lambda$.
\end{theorem}
Matsuzoe \cite{matsuzoe1998realization} proved that if an $n$-dimensional statistical manifold $(\mathbf{M},\nabla,\nabla^{*},g)$ is realized in the affine space $\mathbf{R}^{n+2}$ then it is conformally-projectively flat. Also obtained the converse for simply connected statistical manifolds. \par In the case of non-deenerate centro-affine, equiaffine immersion of codimension two into dually flat statistical manifold we show that $(\mathbf{M},\nabla,h)$ is  conformally-projctively flat.

\begin{theorem}
Let $(\tilde{\mathbf{M}},\tilde{\nabla},\tilde{\nabla}^{*},\tilde{g})$  be a dually flat connected statistical manifold of dimension $(n+2)$ and $\mathbf{M}$ be a manifold of dimension $n$. Let $f: \mathbf{M}\longrightarrow \tilde{\mathbf{M}}$ be a non-degenerate centro-affine, equiaffine immersion of codimension two. Then the inherited statistical manifold $(\mathbf{M},\nabla,h)$ is conformally-projectively flat.
\end{theorem}
\begin{proof}
Let $f: \mathbf{M}\longrightarrow \tilde{\mathbf{M}}$ be a non-degenerate, centro-affine equiaffine immersion of codimension two. For $ p \in \mathbf{M}$ the affine co-ordinate of $f(p)$ in $\tilde{\mathbf{M}}$ with respect to $\tilde{\nabla}$ is denoted by $\theta = [\theta^{i}]$. \par
Let $\eta$ = $\sum_{i=1}^{n+2}\theta^{i}\frac{\partial}{\partial \theta_{i}}$  be the radial vector field of $f$ in $\tilde{\mathbf{M}}$, $\xi$ be the normal vector field of $f$ in $\tilde{\mathbf{M}}$. Then
\begin{eqnarray*}
T_{f(p)}(\tilde{\mathbf{M}}) = f_{*}(T_{p}(\mathbf{M})) + span\lbrace \xi_{p} \rbrace + span\lbrace \eta_{p}\rbrace.
\end{eqnarray*}
Also
\begin{eqnarray*}
 \tilde{\nabla}_{X}f_{*}(Y) &=&  f_{*}(\nabla_{X}Y)+h(X,Y)\xi + T(X,Y)\eta \\
 \tilde{\nabla}_{X}\xi      &=&  -f_{*}(SX)+\tau(X)\xi + \sigma(X)\eta.
\end{eqnarray*}
For a positive function $\psi: \mathbf{M}\longrightarrow \mathbf{R}$, define $g(p) = \psi(p) \theta(f(p))$. Then $(g,\xi,\eta):\mathbf{M}\longrightarrow \tilde{\mathbf{M}}$ is a centro - affine immersion of codimension two. Then \\
 \begin{eqnarray*}
T_{g(p)}(\tilde{\mathbf{M}}) = g_{*}(T_{p}(\mathbf{M})) + span\lbrace \xi_{p} \rbrace + span\lbrace \eta_{p}\rbrace.
\end{eqnarray*}
Also
\begin{eqnarray*}
 \tilde{\nabla}_{X}g_{*}(Y) &=&  g_{*}(\overline{\nabla}_{X}Y)+\overline{h}(X,Y) \xi +\overline{T} (X,Y)\eta \\
 \tilde{\nabla}_{X}\xi      &=&  -g_{*}(\overline{S}X)+
 \overline{\tau}(X)\xi + \overline{\sigma}(X)\eta.
\end{eqnarray*}
Then $(\nabla, h,\tau )$ and $(\overline{\nabla}, \overline{h}, \overline{\tau})$ are related by 
\begin{eqnarray}
\overline{\nabla}_{X}Y  &=& \nabla_{X}Y + d(log \psi)(Y)X+ d(log\psi)(X)Y.\\
\overline{h} &=& \psi h\\
\overline{\tau} &=& \tau. 
\end{eqnarray}
 For each $p\in \mathbf{M}$, choose $\psi$ in some neighbourhood $\mathbf{U}_{p}$ of $p$ such that the coefficient of $\eta$ for the immersion $g=\psi f$ is zero. Take a transversal vector field $\xi^{'}_{x}$ which is equal to $\xi_{p}$ everywhere on $\mathbf{U}_{p}$. Since $g_{*}(T_{x}\mathbf{M}), \eta_{g(x)}$ and $\xi_{x}$ are linearly independent there exist positive function $\phi$, a function $a$ and a tangent vector field $V$ on $\mathbf{U}_{p}$ such that $\phi \xi^{'}_{x} = \xi_{x}+a \eta_{g(x)}+g_{*}V$. \par Now consider the immersion $(g,\xi^{'},\eta): \mathbf{U}_{p}\longrightarrow \tilde{\mathbf{M}}$ then $T^{'} = 0$.\\
Also  
\begin{eqnarray}
\nabla^{'}_{X}Y &=& \overline{\nabla}_{X}Y - \overline{h}(X,Y)V\\
h^{'} &=& \phi \overline{h}\\
\tau^{'}(X)&=& \overline{\tau}(X)- X(log \phi)+ \overline{h}(X,V).
\end{eqnarray}
Since $\xi^{'}$ is parallel around $p$, $\tilde{\nabla}_{X}\xi^{'} = 0$. Then $\mu^{'} = S^{'} = \tau^{'}=0$. Then by Gauss equation $\nabla^{'}$ is flat and so $(\mathbf{U}_{p},\nabla^{'},h^{'})$ is a flat statistical manifold. \par Since $\xi_{x}$ and $\xi^{'}_{x}$ are equiaffine from $(22)$ we have 
\begin{eqnarray*}
\overline{h}(X,V) = X(log \phi).
\end{eqnarray*}
Then 
\begin{equation}
V = \frac{1}{\phi \psi} grad_{h}\phi
\end{equation}
Then from equations $(17)$ to $(23)$ we have, 
\begin{eqnarray*}
h^{'} &=& \phi \psi h\\
h(\nabla^{'}_{X}Y,Z) &=& h(\nabla_{X}Y,Z) - d(log\phi)(Z)h(X,Y) +d(log \psi)(X)h(Y,Z)\\ &+& d(log \psi)(Y)h(X,Z)
\end{eqnarray*}
Hence $(\mathbf{M},\nabla, h)$ is conformally projectively flat.
\end{proof}


\bibliographystyle{spbasic}


\begin{thebibliography}{6}
\providecommand{\natexlab}[1]{#1}
\providecommand{\url}[1]{{#1}}
\providecommand{\urlprefix}{URL }
\expandafter\ifx\csname urlstyle\endcsname\relax
  \providecommand{\doi}[1]{DOI~\discretionary{}{}{}#1}\else
  \providecommand{\doi}{DOI~\discretionary{}{}{}\begingroup
  \urlstyle{rm}\Url}\fi
\providecommand{\eprint}[2][]{\url{#2}}

\bibitem[{Kurose(1990)}]{kurose1990dual}
Kurose T (1990) Dual connections and affine geometry. Mathematische Zeitschrift
  203(1):115--121

\bibitem[{Matsuzoe(1998)}]{matsuzoe1998realization}
Matsuzoe H (1998) On realization of conformally-projectively flat statistical
  manifolds and the divergences. Hokkaido mathematical journal 27(2):409--421

\bibitem[{Amari(2016)}]{amari2016information}
Amari SI (2016) Information Geometry and Its Applications. Springer

\bibitem[{Amari and Nagaoka(2000)}]{amari2000method}
Amari SI, Nagaoka H (2000) Method of Information Geometry. AMS Monograph.
  Oxford University Press Oxford

\bibitem[{Nomizu and Pinkall(1987)}]{nomizu1987geometry}
Nomizu K, Pinkall U (1987) On the geometry of affine immersions. Mathematische
  Zeitschrift 195(2):165--178

\bibitem[{Vos(1989)}]{vos1989fundamental}
Vos PW (1989) Fundamental equations for statistical submanifolds with
  applications to the \uppercase{B}artlett correction. Annals of the Institute
  of Statistical Mathematics 41(3):429--450

\end{thebibliography}

\end{document}